\documentclass{article}

\usepackage{amssymb,amsmath}

\newtheorem{theorem}{Theorem}[section]
\newtheorem{lemma}[theorem]{Lemma}
\newtheorem{fact}[theorem]{Fact}
\newtheorem{corollary}[theorem]{Corollary}

\newcommand{\iin}{\!\in\!}
\newcommand{\rmd}{\,\mathsf{d}}
\newcommand{\abs}[1]{\lvert#1\rvert}
\newcommand{\abslarge}[1]{\left\lvert#1\right\rvert}
\newcommand{\norm}[1]{\lVert#1\rVert}

\newcommand{\ru}{r\,}
\newcommand{\conv}{\star}
\newcommand{\LUC}{\mathsf{LUC}}
\newcommand{\mmeas}{\mathfrak{m}}
\newcommand{\nmeas}{\mathfrak{n}}
\newcommand{\MeasSymb}{{\mathbf{M}}}
\newcommand{\tMeas}{\mathsf{\MeasSymb_t}}
\newcommand{\UMeas}{\mathsf{\MeasSymb_u}}
\newcommand{\BLip}{\mathsf{BLip_b}}
\newcommand{\rtrans}[2]{\rho^{#1}(#2)}

\newcommand{\orbit}[1]{\mathsf{orb}(#1)}
\newcommand{\clorbit}[1]{\overline{\mathsf{orb}}(#1)}
\newcommand{\sect}[1]{\setminus_{#1}}
\newcommand{\tcentre}{\Lambda}
\newcommand{\CBPcentre}{\Lambda^\mathsf{CBP}}
\newcommand{\UMcentre}{\Lambda^\mathsf{URM}}

\newcommand{\mmap}{\varphi}
\newcommand{\mmmap}{\psi}

\newcommand{\expset}{\mathcal{P}}

\newcommand{\psm}{\Delta}
\newcommand{\rp}{\mathsf{RP}}
\newcommand{\ellone}{\ell_1}
\newcommand{\ellinfty}{\ell_\infty}
\newcommand{\wstar}{weak\raisebox{1mm}{$\ast$}\ }
\newcommand{\real}{{\mathbb R}}
\newcommand{\intZ}{{\mathbb Z}}
\newcommand{\qed}{\hfill $\Box$\vspace{3mm}}

\title{Measurable centres in convolution semigroups}

\author{Jan Pachl\thanks{Written while the author was a visitor at the Fields Institute.}
        \\[10pt]
\emph{Fields Institute} \\
\emph{Toronto, Ontario, Canada}}

\date{July 19, 2011}

\begin{document}
\maketitle

\begin{abstract}
In a convolution semigroup over a locally compact group,
measurability of the translation by a fixed element implies continuity.
In other words, the measurable centre coincides with the topological centre.
\end{abstract}


\section{Overview}
    \label{s:overview}

For a topological group $G$, let $\ru G$ denote $G$ with its right uniformity,
and $\LUC(G)$ the space of bounded uniformly continuous functions on $\ru G$.
The norm dual $\LUC(G)^\ast$ of $\LUC(G)$ with convolution $\conv$
is a Banach algebra of some importance in harmonic analysis.
A useful tool for investigating the structure of $\LUC(G)^\ast$ is
its \emph{topological centre}
\[
\tcentre (\LUC(G)^\ast) =  \{ \mmeas \iin \LUC(G)^\ast
\mid \text{ the mapping } \nmeas \mapsto \mmeas \conv \nmeas
\text{ is \wstar continuous on } \LUC(G)^\ast \} .
\]

The space $\tMeas(G)$ of (finite signed) Radon measures on $G$
naturally embeds in $\LUC(G)^\ast$:
A measure $\mu\iin\tMeas(G)$ maps to the functional
$f\mapsto \int f \rmd \mu$, $f\iin\LUC(G)$.
By the theorem of Lau~\cite{Lau1986cam},
$\tcentre (\LUC(G)^\ast)=\tMeas(G)$ for every locally compact group $G$.

In this paper I prove a stronger version of Lau's result, in which \wstar continuity
in the definition of $\tcentre (\LUC(G)^\ast)$ is replaced by measurability.
I also prove similar characterizations for such generalized (measurable) centres
of subsemigroups of $\LUC(G)^\ast$.
In particular, the result applies to the semigroup $\beta G$ for any discrete group $G$,
thus extending a recent result of Glasner~\cite[Th.2.1]{Glasner2009otp}.


\section{Preliminaries}
    \label{s:prel}

All topological spaces and groups considered in this paper are assumed to be Hausdorff,
and all linear spaces to be over the field $\real$ of reals.
Functions (including linear functionals) are real-valued.
It is a simple exercise to extend the results that follow
to linear spaces over the complex field and complex-valued functions.

A pseudometric $\psm$ on a group $G$ is \emph{right-invariant}
iff $\psm(x,y)=\psm(xz,yz)$ for all $x,y$ and $z$ in $G$.
The \emph{right uniformity} on a topological group $G$
is induced by the set of all right-invariant continuous pseudometrics,
denoted by $\rp(G)$.
A bounded pseudometric $\psm$ on $G$ is uniformly continuous in $\ru G$
if and only if there exists $\psm^\prime\iin\rp(G)$ such that $\psm\leq\psm^\prime$.

When $\psm$ is a pseudometric on $G$, write
\[
\BLip(\psm)  =  \{ f\colon G \to \real \mid -1 \leq f(x) \leq 1
\text{ and }
\abs{f(x) - f(y)} \leq \psm(x,y) \text{ for all } x,y \iin G  \} .
\]
Then $\BLip(\psm)$ is a compact subset of the product space $ \real^G $;
in the sequel $\BLip(\psm)$ is always considered with this compact topology.

Let $G$ be a group, $f$ a real-valued function on $G$ and $ x \iin G $.
Define $\rtrans{x}{f}$ (the \emph{right translation} of $f$ by $x$) to be the function
$z\mapsto f(zx)$, $z\iin G$.
The set $\orbit{f} := \{ \rtrans{x}{f} \mid x\iin G \}$ is the \emph{(right) orbit} of $f$.
The closure of $\orbit{f}$ in the product space $\real^G$ is denoted $\clorbit{f}$.
For every $f\iin\LUC(G)$ the set $\orbit{f}$ is norm bounded and uniformly equicontinuous
on $\ru G$, and thus $\clorbit{f}$ is a $G$-pointwise compact subset of $\LUC(G)$.

\begin{fact}[Cor.~15 in~\cite{Pachl2009atg}]
    \label{fact:LCambitable}
Let $G$ be a locally compact group that is not compact.
For every $\psm\iin\rp(G)$ there is $f\iin\LUC(G)$
such that $\BLip(\psm)\subseteq \clorbit{f}$.
\end{fact}

\begin{fact}[\cite{Berezanskii1968mus} and~\cite{Fedorova1967lfd}]
    \label{fact:UMLCGroups}
Let $G$ be a locally compact group and $\mmeas\iin\LUC(G)^\ast$.
If the restriction of $\mmeas$ to the set $\BLip(\psm)$
is continuous for every $\psm\iin\rp(G)$ then $\mmeas\iin\tMeas(G)$.
\end{fact}

By combining the first two facts we obtain a characterization of finite Radon measures
on locally compact groups:
A functional $\mmeas\iin\LUC(G)^\ast$ is in $\tMeas(G)$ if and only if
for every $f\iin\LUC(G)$ the restriction of $\mmeas$ to $\clorbit{f}$ is
$G$-pointwise continuous.

Let $X$ be a (Hausdorff as always) topological space and $A\subseteq X$.
Say that $A$ is a \emph{CBP set} iff
for every continuous mapping $\mmap\colon K\to X$
from a compact space $K$ the set $\mmap^{-1}(A)$ has the Baire property in $K$.
Say that a real-valued function $f$ on $X$ is \emph{CBP measurable}
iff $f^{-1}(U)$ is a CBP set in $X$ for every open subset $U$ of $\real$.

When $X$ is compact, CBP subsets of $X$ are exactly the universally Baire-property sets
in the terminology of Fremlin~\cite{Fremlin2011tsf}.
Evidently the CBP subsets of $X$ form a $\sigma$-algebra,
every Borel set is CBP, and every Borel measurable mapping is CBP measurable.
If a mapping $f\colon X \to Y$ is CBP measurable then so is its restriction to
any subspace of $X$.

Denote by $\expset I$ the set of all subsets of a set $I$, and identify
$\expset I$ with the compact set $2^I$.

\begin{fact}[Lemma~2.1 in~\cite{Christensen1981mff}]
    \label{fact:BPmeasure}
Let $I$ be an infinite set,
and let $\mu\colon \expset I \to \real$ be finitely additive.
If $\mu$ is  CBP measurable then it is a measure on $\expset I$
and $\mu(A)=\sum_{i\in A} \mu(\{i\})$ for every $A\subseteq I$.
\end{fact}

The next theorem and its proof are due to Fremlin.
A slightly different version appears in~\cite[1E]{Fremlin2011tsf}.

\begin{theorem}
    \label{th:CBPquotient}
Let $K_0$ be a compact space and $\mmap_0\colon K_0 \to X$ a continuous surjective mapping
onto a compact space $X$.
Let $A\subseteq X$ be such that $\mmap_0^{-1}(A)$ is a CBP set in~$K_0$.
Then $A$ is a~CBP set in $X$.
\end{theorem}

\noindent
\textbf{Proof}.
Take any continuous mapping $\mmap_1\colon K_1 \to X$ from a compact space $K_1$.
For $i=0,1$ let $\pi_i \colon K_0 \times K_1 \to K_i$ be the canonical projections.
Then
\[
K:=\{ (x_0,x_1 ) \iin K_0 \times K_1 \mid \mmap_0(x_0) = \mmap_1(x_1) \}
\]
is a compact subset of $K_0 \times K_1$,
and $\pi_1(K)=K_1$ because $\mmap_0(K_0)=X$.
By 4A2Gi in~\cite{Fremlin2006mt4}
there is a closed set $K^\prime \subseteq K$
such that the restriction of $\pi_1$ to $K^\prime$ is irreducible and
$\pi_1(K^\prime)=K_1$.

By ~\cite[L.2]{Pol1976rrb}, \cite[25.2.3]{Semadeni1971bsc},
if a set has the Baire property in a compact space,
then so does its image under any irreducible continuous surjection.
Since $\mmap_0^{-1}(A)$ is a CBP set,
the set $K^\prime\cap \pi_0^{-1} (\mmap_0^{-1}(A))$ has the Baire property in $K^\prime$,
and the set
\[
\pi_1(K^\prime\cap \pi_0^{-1} (\mmap_0^{-1}(A))) = \mmap_1^{-1} (A)
\]
has the Baire property in $K_1$.
\qed

The following theorem is an essential step in the proof of the main result in the next section.
As before, $\BLip(\psm)$ is considered with the $G$-pointwise topology.

\begin{theorem}
    \label{th:MeasMeasure}
Let $G$ be a locally compact group and $\mmeas\iin\LUC(G)^\ast$,
and assume that for every $\psm\iin\rp(G)$
the restriction of $\mmeas$ to $\BLip(\psm)$ is CBP measurable.
Then for every $\psm\iin\rp(G)$ the restriction of $\mmeas$ to $\BLip(\psm)$ is continuous.
\end{theorem}

For metrizable locally compact groups,
Theorem~\ref{th:MeasMeasure} is a direct consequence of Theorem~2 and Lemma~4.1
in~\cite{Christensen1981mff}.
Here I prove the general case, after several auxiliary lemmas.

A \emph{partition of unity} on a set $X$ is a mapping $p\colon X\to \ellinfty(I)$
where $I$ is a non-empty index set,
$0\leq p(x)(i)\leq 1$ for all $x\iin X$ and $i\iin I$,
and $\sum_{i\in I} p(x)(i) = 1 $ for every $x\iin X$.
Write $p_i(x) := p(x)(i)$;
thus each $p_i$ is a function on $X$ with values in the interval $[0,1]$.
Note that the range of $p$ is included in $\ellone(I)\subseteq\ellinfty(I)$.
Denote  the $\ellone(I)$ norm by $\norm{\cdot}_1$ and
the $\ellinfty(I)$ norm by $\norm{\cdot}_\infty$.

When $\psm$ is a pseudometric on $X$,
say that the partition of unity $p\colon X\to \ellinfty(I)$
is \emph{subordinated to $\psm$} iff for every $i\iin I$ we have $\psm(x,y)\leq 1$
whenever $x,y\iin G$, $p_i(x)>0$, $p_i(y)>0$.
When $G$ is a topological group, say that the uniform space $\ru G$ has the
\emph{$(\ellone)$ property} iff for every $\psm\iin\rp(G)$
there exists a partition of unity $p$ on $G$
that is subordinated to $\psm$ and
uniformly continuous from $\ru G$ to $\ellone(I)$ with the $\norm{\cdot}_1$ norm.

\begin{lemma}
    \label{lemma:elloneprop}
For every locally compact group $G$ the uniform space $\ru G$ has the ($\ellone$) property.
\end{lemma}

For metrizable locally compact groups this is Lemma~4.1 in~\cite{Christensen1981mff}.
Essentially the same proof works for the general case, and I do not repeat it here. \\

\begin{lemma}
    \label{lemma:MeasSum}
Let $G$ be a topological group and let $\mmeas\iin\LUC(G)^\ast$ be
such that for every $\psm\iin\rp(G)$ the restriction of $\mmeas$ to $\BLip(\psm)$
is CBP measurable.
Let $I$ be a non-empty index set and $\mmap\colon \ru G \to \ellone(I)$
a uniformly continuous mapping from $\ru G$ to $\ellone(I)$ with the $\norm{\cdot}_1$ norm,
and such that $\norm{\mmap(x)}_1 \leq 1$ for every $x\iin G$.
Then
\[
\sum_{i\in I} \,\abs{\mmeas (\mmap_i)} < \infty \;\;\;\text{ and }\;\;\;
\sum_{i\in A} \mmeas (\mmap_i)
= \mmeas \left( \sum_{i\in A} \mmap_i \right)
\;\;\text{ for }\;\; A\subseteq I
\]
where $\mmap_i(x):=\mmap(x)(i)$ for $i\iin I$, $x\iin X$.
\end{lemma}

\noindent
\textbf{Proof}.
Since $\mmap$ is uniformly continuous in the $\norm{\cdot}_1$ norm,
there is $\psm\iin\rp(G)$ such that $\norm{\mmap(x)-\mmap(y)}_1 \leq \psm(x,y)$
for all $x,y\iin G$.
The expression
\[
\mmmap(A) := \sum_{i\in A} \mmap_i(x), \;\; A\subseteq I,
\]
defines a continuous finitely additive mapping $\mmmap\colon \expset I \to \BLip(\psm) $.
Hence the function $\mmeas\circ\mmmap$ is finitely additive and CBP measurable on $\expset I$.
Apply Fact~\ref{fact:BPmeasure}.
\qed

\noindent
\textbf{Proof of Theorem~\ref{th:MeasMeasure}}.
Take any $\psm\iin\rp(G)$ and any net $\{f_\gamma\}_\gamma$
of functions $f_\gamma \iin \BLip(\psm)$ such that $\lim_\gamma f_\gamma(x) = 0$
for all $x\iin G$.
Fix an arbitrary $\varepsilon>0$.
By Lemma~\ref{lemma:elloneprop} there is a partition of unity $p$ on $G$
that is subordinated to $\psm/\varepsilon$ and
uniformly continuous from $\ru G$ to $\ellone(I)$ with the $\norm{\cdot}_1$ norm.

For each $i\iin I$ choose a point $x_i\iin G$ such that $\psm(x,x_i)\leq\varepsilon$
whenever $p_i(x)> 0$.
Then
\[
\abslarge{f_\gamma(x) - \sum_{i\in I} f_\gamma (x_i) \cdot p_i (x)}
\leq \sum_{i\in I} p_i (x) \cdot \abs{f_\gamma(x) -  f_\gamma (x_i)}
\leq \sum_{i\in I} p_i (x) \cdot \psm(x,x_i)
\leq \varepsilon
\]
for all $\gamma$ and all $x\iin G$.

For a fixed $\gamma$, define $\mmap\colon G \to \ellone(I)$ by
$\mmap(x)(i):=f_\gamma(x_i)\cdot p_i(x)$, $x\iin X$, $i\iin I$,
and apply Lemma~\ref{lemma:MeasSum} to get
\[
\sum_{i\in I} f_\gamma(x_i) \mmeas (p_i)
= \mmeas \left( \sum_{i\in I} f_\gamma(x_i) p_i \right) .
\]

By Lemma~\ref{lemma:MeasSum} there is a finite set ${D}\subseteq I$ such that
$\sum_{i\in I\setminus {D}} \,\abs{\mmeas (p_i)} < \varepsilon$.
For almost all $\gamma$ we have $\abs{f_\gamma (x_i)} < \varepsilon$ when $i\iin {D}$, and
\begin{align*}
\abs{\mmeas(f_\gamma)}
    & \leq \abslarge{\,\mmeas\left(f_\gamma
                        - \sum_{i\in I} f_\gamma (x_i) \cdot p_i \right)}
        + \abslarge{\,\sum_{i\in I} f_\gamma(x_i) \mmeas(p_i)\,} \\
    & \leq \norm{\mmeas} \varepsilon
        + \sum_{i\in {D}} \,\abs{f_\gamma (x_i)} \cdot \abs{\mmeas(p_i)}
        + \sum_{i\in I\setminus {D}} \abs{f_\gamma (x_i)} \cdot \abs{\mmeas(p_i)}   \\
    & \leq \norm{\mmeas} \varepsilon
        + 2 \norm{\mmeas} \varepsilon + \varepsilon
        = (3 \norm{\mmeas} + 1 ) \varepsilon  .
\end{align*}
As this holds for every $\varepsilon>0$, $\mmeas$ is continuous on $\BLip(\psm)$.
\qed


\section{Generalized centres}
    \label{s:centres}

For any topological group $G$, the convolution in $\LUC(G)^\ast$ may be written as
\[
\mmeas \conv \nmeas (f)  = \mmeas ( \sect{x} \nmeas ( \sect{y} f(xy ) ) )
\]
for $\mmeas,\nmeas\iin\LUC(G)^\ast$ and $f\iin\LUC(G)$.
Here $\sect{x} f(\dotsc)$ means ``$f(\dotsc)$ as a function of $x$''.
This formula applies not only in $\LUC(G)^\ast$ for a topological group $G$
but also in analogous spaces over more general
\emph{semiuniform semigroups}~\cite{Pachl2008ssc}.

$\LUC(G)^\ast$ with convolution is a Banach algebra.
Here we mostly treat it as a semigroup with the $\conv$ operation.
The group $G$ naturally embeds in $\LUC(G)^\ast$:
An element $x\iin G$ maps to the functional $f\mapsto f(x)$, $f\iin\LUC(G)$.
The embedding is a homeomorphism of $G$ onto its image in $\LUC(G)^\ast$ with the \wstar
topology.
The embedding also preserves the algebraic structure, so that $G$ may be identified
with a subgroup of $\LUC(G)^\ast$.

The \wstar closure of $G$ in $\LUC(G)^\ast$, denoted here $G^\LUC$,
is a \wstar compact subsemigroup of $\LUC(G)^\ast$.
It is known as
the \emph{canonical $\mathcal{LC}$-compactification}~\cite{Berglund1989aos},
\emph{universal enveloping semigroup}~\cite{deVries1993etd},
\emph{$\mathcal{LUC}$-compactification}~\cite{Lau1995tcc},
or \emph{greatest ambit}~\cite{Pestov2006did} of $G$;
or, in the language of uniform spaces,
a \emph{uniform (or Samuel) compactification} of $\ru G$.
When $G$ is discrete, $G^\LUC$ is its \v{C}ech--Stone compactification $\beta G$.
When $G$ is locally compact, $G^\LUC \cap \tMeas(G) = G$.

For any topological group $G$ and any $S \subseteq \LUC(G)^\ast$ define
\begin{align*}
\tcentre (S) := \{ \mmeas \iin S \mid \forall f\iin\LUC(G)
               & \text{ the function } \nmeas \mapsto \mmeas \conv \nmeas (f)   \\
               & \text{ is \wstar continuous on } S \}                     \\
\CBPcentre (S) := \{ \mmeas \iin S \mid \forall f\iin\LUC(G)
               & \text{ the function } \nmeas \mapsto \mmeas \conv \nmeas (f)   \\
               & \text{ is \wstar CBP measurable on } S \}
\end{align*}
(the \emph{topological centre} and the \emph{\wstar CBP measurable centre of $S$}).

It is well known and easy to prove that
$S\cap\tMeas(G)\subseteq\tcentre(S)\subseteq\CBPcentre(S)\subseteq S$.
If $G$ is compact then $\LUC(G)^\ast=\tMeas(G)$
and therefore $\tcentre(S)=\CBPcentre(S)=S=S\cap\tMeas(G)$
for every $S\subseteq\LUC(G)^\ast$.

Now we come to the main result of this paper.
The proof strategy is the same as in section~5 of~\cite{Pachl2009atg}.

\begin{theorem}
    \label{th:UBPcentre}
Let $G$ be a locally compact group and $G^\LUC\subseteq S \subseteq \LUC(G)^\ast$.
Then
\[
\tcentre(S)=\CBPcentre(S)=S\cap\tMeas(G).
\]
\end{theorem}

\noindent
\textbf{Proof}.
In view of the preceding discussion,
it is enough to prove that $\CBPcentre(S)\subseteq S\cap\tMeas(G)$ when
$G$ is not compact.

For $f\iin\LUC(G)$ define the mapping $\mmap_f\colon\LUC(G)^\ast \to \LUC(G)$ by
\[
\mmap_f(\nmeas):= \sect{x} \nmeas (\sect{y} f(xy)), \;\;\nmeas\iin\LUC(G)^\ast.
\]
Then for every $\mmeas\iin\LUC(G)^\ast$ the mapping $\nmeas \mapsto \mmeas \conv \nmeas (f)$
from $\LUC(G)^\ast$ to $\real$ is the composition $\mmeas \circ \mmap_f$.
By~\cite[L.19]{Pachl2009atg}, $\mmap_f$ is continuous from $G^\LUC$ to $\LUC(G)$ with
the $G$-pointwise topology, and $\mmap_f(G^\LUC)=\clorbit{f}$.

\begin{picture}(300,90)(-80,0)
\thicklines
\put(50,60){\makebox(30,20){$G^\LUC$}}
\put(140,60){\makebox(40,20){$\clorbit{f}$}}
\put(140,0){\makebox(40,20){$\real$}}

\put(80,70){\vector(1,0){60}}
\put(160,60){\vector(0,-1){40}}
\put(80,60){\vector(3,-2){66}}

\put(95,66){\makebox(30,20){$\mmap_f$}}
\put(160,30){\makebox(20,20){$\mmeas$}}
\put(40,20){\makebox(80,20){$\nmeas \mapsto \mmeas \conv \nmeas (f) $}}
\end{picture}

Now assume that $\mmeas\iin\CBPcentre(S)$,
which means that for every $f\iin\LUC(G)$ the mapping $\mmeas \circ \mmap_f$
is CBP measurable on $S$, and therefore also on $G^\LUC \subseteq S$.
By Theorem~\ref{th:CBPquotient},
$\mmeas$ is CBP measurable on $\clorbit{f}$.
By Fact~\ref{fact:LCambitable},
$\mmeas$ is  CBP measurable on $\BLip(\psm)$
for every right-invariant continuous pseudometric $\psm$ on $G$,
and therefore also continuous on $\BLip(\psm)$ by Theorem~\ref{th:MeasMeasure}.
Hence $\mmeas\iin\tMeas(G)$ by Fact~\ref{fact:UMLCGroups}.
\qed

By choosing $S=\LUC(G)^\ast$ and $S=G^\LUC$ we obtain two corollaries.
The first one is the promised strengthening of Lau's theorem~\cite{Lau1986cam}.

\begin{corollary}
    \label{cor:tcentredual}
$\tcentre(\LUC(G)^\ast)=\CBPcentre(\LUC(G)^\ast)=\tMeas(G)$
for every locally compact group $G$.
\end{corollary}

The second corollary is a common generalization
of the theorems of Lau and Pym~\cite{Lau1995tcc}
and Glasner~\cite{Glasner2009otp}.

\begin{corollary}
    \label{cor:tcencomp}
$\tcentre(G^\LUC)=\CBPcentre(G^\LUC)=G$
for every locally compact group $G$.
\end{corollary}

Note that Theorem~\ref{th:UBPcentre} applies also to many other
sets between $G^\LUC$ and $\LUC(G)^\ast$
--- for example, the set of positive elements in $\LUC(G)^\ast$,
or the set of means on $\LUC(G)$.


\section{Variations and open problems}
    \label{s:extensions}

One may ask to what extent Theorem~\ref{th:UBPcentre} and its corollaries depend
on the group $G$ being locally compact.
The results in~\cite{Ferri2007tca} and~\cite{Pachl2009atg} suggest that the space
$\UMeas(\ru G)$ of \emph{uniform measures} should take the place
of $\tMeas(G)$ in describing the centres in convolution semigroups as we move beyond
locally compact groups ($\UMeas(\ru G)$ and $\tMeas(G)$ coincide when $G$ is locally compact).
This leads to the question whether $\CBPcentre(\LUC(G)^\ast)=\UMeas(\ru G)$
for every topological group $G$,
or at least for some interesting class of non-locally-compact groups.

With the same approach as in the proof of Theorem~\ref{th:UBPcentre},
we get that $\CBPcentre(S)=S\cap\UMeas(\ru G)$ for $G^\LUC\subseteq S \subseteq \LUC(G)^\ast$
whenever $G$ is an ambitable topological group~\cite{Pachl2009atg} for which
$\ru G$ has the $(\ellone)$ property.
However,
infinite-dimensional normed spaces do not have the ($\ellone$) property
by the theorem of Zahradn{\'{\i}}k~\cite{Zahradnik1976lcp}.

One may also try to weaken the measurability condition in the definition of
$\CBPcentre(S)$.
In one direction, Schachermayer's example~\cite{Schachermayer1981mcl}
marks a limit of such generalizations:
For the additive group $c_0$ and the metric $\psm$ of the sup norm on $c_0$,
there is a bounded linear functional $\mmeas$ on $\LUC(c_0)$ whose restriction
to $\BLip(\psm)$ is Baire-property measurable and yet $\mmeas$ is not in $\tMeas(c_0)$.

In another direction, for $S \subseteq \LUC(G)^\ast$ define
\begin{align*}
\UMcentre (S) = \{ \mmeas \iin S \mid \forall f\iin\LUC(G)
               & \text{ the function } \nmeas \mapsto \mmeas \conv \nmeas (f)   \\
               & \text{ is \wstar universally Radon-measurable on } S \} .
\end{align*}

The characterization of $\UMcentre(S)$ is not as straightforward as that of $\CBPcentre(S)$,
even for the group $\intZ$ of integers with the discrete topology.
On one hand, Glasner's proof of Theorem~2.1 in~\cite{Glasner2009otp} demonstrates that
if $G$ is a countable discrete group then $\UMcentre(\beta G)=G$,
which improves (for such groups) Corollary~\ref{cor:tcencomp}.
On the other hand,
the statement $\UMcentre(\LUC(\intZ)^\ast)=\tMeas(\intZ)$,
which may be written simply as $\UMcentre(\ellinfty^\ast)=\ellone$,
is neither provable nor disprovable in the ZFC set theory.
That follows from old results about \emph{medial limits},
covered by Fremlin~\cite[538Q]{Fremlin2008mt5},
along with a recent result of Larson~\cite{Larson2009fdm}.

\vspace{5mm}
\textbf{Acknowledgements}.
I wish to thank David Fremlin for the proof of Theorem~\ref{th:CBPquotient},
and Petr Holick\'{y} for relevant references.
The questions addressed in this paper originated in discussions with Matthias Neufang,
Dilip Raghavan and Juris Stepr\={a}ns.
I appreciate their input.



\begin{thebibliography}{10}

\bibitem{Berezanskii1968mus}
Berezanski{\u\i}, I.~A.
\newblock \textit{ Measures on uniform spaces and molecular measures\/}.
\newblock Trudy Moskov. Mat. Ob\v s\v c. \textbf{19} (1968), 3--40.
\newblock English translation: Trans. Moscow Math. Soc. 19 (1968), 1--40.

\bibitem{Berglund1989aos}
Berglund, J.~F., Junghenn, H.~D., and Milnes, P.
\newblock \textit{ Analysis on semigroups}.
\newblock John Wiley \& Sons Inc., New York, 1989.

\bibitem{Christensen1981mff}
Christensen, J. P.~R., and Pachl, J.
\newblock \textit{ Measurable functionals on function spaces\/}.
\newblock Ann. Inst. Fourier (Grenoble) \textbf{31}, 2 (1981), 137--152.

\bibitem{deVries1993etd}
de~Vries, J.
\newblock \textit{ Elements of topological dynamics}, Vol.~257.
\newblock Kluwer Academic Publishers Group, Dordrecht, 1993.

\bibitem{Fedorova1967lfd}
Fedorova, V.~P.
\newblock \textit{ Linear functionals and {D}aniell integral on spaces of
  uniformly continuous functions\/}.
\newblock Mat. Sb. (N.S.) \textbf{74 (116)} (1967), 191--201.
\newblock English translation: Math. USSR -- Sbornik 3 (1967), 177--185.

\bibitem{Ferri2007tca}
Ferri, S., and Neufang, M.
\newblock \textit{ On the topological centre of the algebra {${\rm LUC}(G)\sp \ast$}
  for general topological groups\/}.
\newblock J. Funct. Anal. \textbf{244}, 1 (2007), 154--171.

\bibitem{Fremlin2006mt4}
Fremlin, D.~H.
\newblock \textit{ Measure theory. {V}ol. 4, Topological measure spaces. Parts I, II}.
\newblock Torres Fremlin, Colchester.
\newblock Corrected second printing, 2006.

\bibitem{Fremlin2008mt5}
Fremlin, D.~H.
\newblock \textit{ Measure Theory. {V}ol. 5, Set-theoretic measure theory. Parts I, II}.
\newblock Torres Fremlin, Colchester.
\newblock 2008.

\bibitem{Fremlin2011tsf}
Fremlin, D.~H.
\newblock \textit{ {Topological spaces after forcing (16.6.11)}\/}.
\newblock \\ http://www.essex.ac.uk/maths/people/fremlin/n05622.ps.

\bibitem{Glasner2009otp}
Glasner, E.
\newblock \textit{ On two problems concerning topological centers\/}.
\newblock Topology Proc. \textbf{33} (2009), 29--39.

\bibitem{Larson2009fdm}
Larson, P.~B.
\newblock \textit{ The filter dichotomy and medial limits\/}.
\newblock J. Math. Log. \textbf{9}, 2 (2009), 159--165.

\bibitem{Lau1986cam}
Lau, A. T.-M.
\newblock \textit{ Continuity of {A}rens multiplication on the dual space of
  bounded uniformly continuous functions on locally compact groups and
  topological semigroups\/}.
\newblock Math. Proc. Cambridge Philos. Soc. \textbf{99}, 2 (1986), 273--283.

\bibitem{Lau1995tcc}
Lau, A. T.-M., and Pym, J.
\newblock \textit{ The topological centre of a compactification of a locally
  compact group\/}.
\newblock Math. Z. \textbf{219}, 4 (1995), 567--579.

\bibitem{Pachl2008ssc}
Pachl, J.
\newblock \textit{ {Semiuniform semigroups and convolution}\/}.
\newblock arXiv (2008), arXiv:0811.3576.

\bibitem{Pachl2009atg}
Pachl, J.
\newblock \textit{ Ambitable topological groups\/}.
\newblock Topology Appl. \textbf{156}, 13 (2009), 2200--2208.

\bibitem{Pestov2006did}
Pestov, V.
\newblock \textit{ Dynamics of infinite-dimensional groups}, University Lecture
  Series, Vol.~40.
\newblock American Mathematical Society, Providence, RI, 2006.

\bibitem{Pol1976rrb}
Pol, R.
\newblock \textit{ Remark on the restricted {B}aire property in compact spaces\/}.
\newblock Bull. Acad. Polon. Sci. S\'er. Sci. Math. Astronom. Phys.
  \textbf{24}, 8 (1976), 599--603.

\bibitem{Schachermayer1981mcl}
Schachermayer, W.
\newblock \textit{ Measurable and continuous linear functionals on spaces of
  uniformly continuous functions\/}.
\newblock Measure theory ({O}berwolfach, 1981), Lecture Notes in Math.,
  Vol.~945. Springer, Berlin, 1982, pp.~155--166.

\bibitem{Semadeni1971bsc}
Semadeni, Z.
\newblock \textit{ Banach spaces of continuous functions. {V}ol. {I}}.
\newblock PWN---Polish Scientific Publishers, Warsaw, 1971.
\newblock Monografie Matematyczne, Tom 55.

\bibitem{Zahradnik1976lcp}
Zahradn{\'{\i}}k, M.
\newblock \textit{ {$l\sb{1}$}-continuous partitions of unity on normed spaces\/}.
\newblock Czechoslovak Math. J. \textbf{26(101)}, 2 (1976), 319--329.

\end{thebibliography}
\end{document}